\documentclass[12pt]{article}
\usepackage{amsmath,amssymb,color,amscd}
\usepackage{xspace}

\newcommand{\Var}{{\rm{Var}_{\mathbb{C}}}}

\newcommand{\ind}{{\mbox{\rm ind}}}

\newcommand{\chiorb}{{\chi^{\rm orb}}}

\def\1{\underline{1}}

\def\Z{{\mathbb Z}}

\def\Q{{\mathbb Q}}

\def\ind{{\rm ind}}

\newtheorem{theorem}{Theorem}

\newtheorem{lemma}{Lemma}
\newtheorem{proposition}{Proposition}
\newtheorem{definition}{Definition}

\newenvironment{remark}
{\smallskip\noindent{\bf Remark\/}.}{\smallskip\par}

\newenvironment{proof}
{\noindent{\bf Proof\/}.}{{ $\square$}\smallskip\par}

\title{Tamanoi equation for orbifold Euler characteristics: revisited
}

\author{S.M.~Gusein-Zade \thanks{The paper is an output of a research project implemented as part of the Basic Research Program at the National Research University Higher School of Economics (HSE University).
Address: Lomonosov Moscow State University, Faculty
of Mechanics and Mathematics and 
Moscow Center for Fundamental and Applied Mathematics,
GSP-1, Moscow, 119991, Russia \&
National Research University ``Higher School of Economics'',
Usacheva street 6, Moscow, 119048, Russia. E-mail:
sabir\symbol{'100}mccme.ru} }
\date{}
\begin{document}
\def\eps{\varepsilon}

\maketitle

\centerline{To Sergei Petrovich with deep gatitute}

\begin{abstract}
Tamanoi equation is a Macdonald type equation for the orbifold Euler characteristic and for its analogues of higher orders. It claims that the generating series of the orbifold Euler characteristics of a fixed order of analogues of the symmetric powers for a space with a finite group action can be represented as a certain unified (explicitly written) power series in the exponent equal to the orbifold Euler characteristic of the same order of the space itself. In the paper,
in particular, we explain how the Tamanoi equation follows from its verification for actions of (finite) groups on the one-point space.
Statements used for that are generalized to analogues of the orbifold
Euler characteristic corresponding to finitely generated groups.
It is shown that, for these generalizatins, the analogue of the Tamanoi equation does not hold in general.
\end{abstract}

\section{Introduction}\label{sec:intro}
In all the paper, $\chi(\cdot)$ denotes {\em the additive Euler characteristic} defined as the alternating sum of the ranks of the cohomology groups with compact support:
$$
\chi(X)=\sum_{q\ge 0}(-1)^q \dim H^q_c(X;\Z)\,.
$$
The well-known Macdonald equation
(\cite{Macdonald}, see also~\cite[Theorem~1]{RMS-2017}) states that,
for a (nice enough) topological space $X$ one has the equality \begin{equation}
 1+\sum_{k=1}^{\infty} \chi(S^kX)\cdot t^k=(1-t)^{-\chi(X)},
\end{equation}
wher $S^kX=X^k/S_k$ is the $k$th symmetric power of the space $X$ ($S_k$ is the permutation group on $k$ elements).

For a space $X$ with an action of a finite group $G$,
there is defined the orbifold Euler characteristic (see, e.\,g., \cite{AS}, \cite{HH}).
Let ${\rm Conj\,}G$ be the set of the conjugacy classes of the elements of the group $G$, for an element $g\in G$ let $C_G(g)$ denote
the centralizer $\{h\in G: h^{-1}gh=g\}$ of $g$. For the $G$-space
$X$ and a subgroup $H\subset G$, $X^H$ denotes the set
$\{x\in X: gx=x \text{\ \ for all\ \ } g\in H\}$ of fixed points of the subgroup $H$.
{\em The orbifold Euler characteristic}
$\chiorb(X,G)$ of a $G$-space $X$ is defined by the (equivalent: see,
e.\,g.,~\cite[Proposition~6]{RMS-2017}) equations:
\begin{equation}\label{eqn:chi-orb}
\chiorb(X,G)=\frac{1}{\vert G\vert}
 \sum_{{(g_0,g_1)\in G^2:}\atop{g_0g_1=g_1g_0}}\chi(X^{\langle g_0,g_1\rangle})
=\sum_{[g]\in {{\rm Conj\,}G}} \chi(X^{\langle g\rangle}/C_G(g))\,,
\end{equation}
where $g$ is a representative of the class $[g]$, $\langle \ldots\rangle$ is the subgroup of $G$ generated by the corresponding elements.

For a $G$-space $X$, its {\em orbifold Euler characteristics of higher orders} where defined in~\cite{AS} and~\cite{Bryan_Fulman-1998} by the equations:
\begin{equation}\label{eqn:chi-k-orb}
 \chi^{(k)}(X,G)=
\frac{1}{\vert G\vert}\sum_{{{\bf g}\in G^{k+1}:}\atop{g_ig_j=g_jg_i}}\chi(X^{\langle {\bf g}\rangle})
=\sum_{[g]\in {{\rm Conj\,}G}} \chi^{(k-1)}(X^{\langle g\rangle}, C_G(g))\,,
\end{equation}
wher $k$ is a non-negative integer (the order of the Euler characteristic), ${\bf g}=(g_0,g_1, \ldots, g_k)$
and (for the second~--- recurrent definition)
$\chi^{(0)}(X,G)$ is given by the first definition
and coincides with the Euler characteristic
$\chi(X/G)$ of the quotient: see, e.\,g., \cite[Proposition~7]{RMS-2017}).
The usual orbifold Euler characteristic $\chiorb(\cdot,\cdot)$
is the orbifold Euler characteristic of order 1. (Into this series
one may include the Euler-Satake characteristic (\cite{Satake}), which can be interpreted in a certain (not completely streightforward) way as  $\chi^{(-1)}(X,G)$.)

One has generalizations of the orbifold Euler characteristics of higher orders corresponding to (arbitrary) finitely generated groups,
for whom the characteristics $\chi^{(k)}$ correspond to free abelian groups:~\cite{Tamanoi}, \cite{Farsi_Seaton-2011}. Let $A$ be
a finitely generated group. {\em The orbifold Euler characteristic correspondingto the group $A$},
(or {\em the $A$-Euler characteristic}) is
\begin{equation}\label{eqn:chi(A)}
\chi^{(A)}(X,G)=
\frac{1}{\vert G\vert}\sum_{\varphi\in {\rm Hom\,}(A,G)}\chi(X^{\langle {\varphi(A)}\rangle})\,,
\end{equation}
where the summation is over the set of homomorphisms
$\varphi:A\to G$, $X^{\langle {\varphi(A)}\rangle}$ is the set
of point of $X$ fixed with respect to all the elements of the image
of $\varphi$. (The value $\chi^{(A)}(X,G)$ is, in general,
a rational number.
If $A=A'\oplus \Z$, them $\chi^{(A)}(X,G)$ is an integer: Proposition~\ref{prop:A-Euler_for_product}, applied to $A_1=\Z$, $A_2=A'$.
One has $\chi^{(k)}(X,G)=\chi^{(\Z^{k+1})}(X,G)$.)
In addition
\begin{equation}
 \chi^{(A)}(G/G,G)=\frac{1}{\vert G\vert}
 \left| {\rm Hom}(A,G)\right|\label{eqn:chi_A_G/G}\,.
\end{equation}

It is not difficult to see that, for the trivial group $G$ (i.\,e., $G=\{e\}$)
the $A$-Euler characteristic $\chi^{(A)}(X,\{e\})$
coincides with the usual Euler characteristic
$\chi(X)$. The definition directly implies
the additivity and the multiplicativity of the $A$-Euler characteristic
in the following sense: if $X_1$ and $X_2$ are $G$-spaces, then
$\chi^{(A)}(X_1\sqcup X_2,G)=\chi^{(A)}(X_1,G)+\chi^{(A)}(X_2,G)$; if $X_1$ is a $G_1$-space and $X_2$is a $G_2$-space, then
\begin{equation}\label{eqn:multiplicat}
\chi^{(A)}(X_1\times X_2,G_1\times G_2)=\chi^{(A)}(X_1,G_1)\cdot\chi^{(A)}(X_2,G_2).
\end{equation}

\begin{proposition}\label{prop:A-Euler_for_product} (\cite[Proposition~2-1]{Tamanoi})
 If $A=A_1\times A_2$, then
 \begin{equation}\label{eqn;A-Euler_for_product}
  \chi^{(A)}(X,G)=\sum_{[\varphi]\in {\rm Hom}(A_1, G)/G}
  \chi^{(A_2)}\left(X^{\langle\varphi(A_1)\rangle},C_G(\varphi(A_1))\right)\,,
 \end{equation}
 where the group $G$ acts on the set ${\rm Hom}(A_1, G)$ by conjugation,
 $\varphi$ is a representative of the conjugacy class $[\varphi]$.
\end{proposition}

\begin{proof} (c.\,f.\ the proofs of \cite[Proposition~2-1]{Tamanoi} and ~\cite[Proposition~6]{RMS-2017}) One has:
\begin{eqnarray*}
 \ &{\ }&\chi^{(A)}(X,G)=\frac{1}{\vert G\vert}
 \sum_{\varphi\in {\rm Hom}(A, G)}\chi(X^{\langle {\varphi(A)}\rangle})=\\
 \ &=&\frac{1}{\vert G\vert}\sum_{[\varphi]\in{\rm Hom}(A_1,G)/G}
 \vert[\varphi]\vert\cdot\sum_{\psi\in{\rm Hom}(A_2,C_G(\varphi(A_1))}
 \chi(X^{\langle\varphi(A_1)\cdot\psi(A_2)\rangle})=\\
 \ &=&\frac{1}{\vert G\vert}\sum_{[\varphi]\in{\rm Hom}(A_1,G)}
 \frac{\vert G\vert}{\vert C_G(\varphi(A_1))\vert}
 \sum_{\psi\in{\rm Hom}(A_2,C_G(\varphi(A_1))}
 \chi((X^{\langle\varphi(A_1)\rangle})^{\langle\psi(A_2)\rangle})=\\
 \ &=&\sum_{[\varphi]\in {\rm Hom}(A_1, G)/G}
  \chi^{(A_2)}(X^{\langle\varphi(A_1)\rangle},C_G(\varphi(A_1)))\,.
\end{eqnarray*}
\end{proof}

As an analogue of the $n$th symmetric power of a space $X$ with an
action of a (finite) group $G$ one can consider the Cartesian power $X^n$ of $X$ with the action of the wreath product
$G_n=G\wr S_n=G^n\rtimes S_n$. 
The multiplication in $G_n$ is given by the formula
$$((g_1,\ldots,g_n),s)((g'_1,\ldots,g'_n),s')=
((g_1g'_{s^{-1}(1)},\ldots g_ng'_{s^{-1}(n)}),ss').$$
The action of the group $G_n$ on the Cartesian power $X^n$ is given by the formula
$$
((g_1,\ldots,g_n),s)(x_1,\ldots, x_n)=
(g_1x_{s^{-1}(1)}, \ldots, g_nx_{s^{-1}(n)})\,,
$$
where $x_i\in X$, $g_i\in G$, $s\in S_n$.
(One has $X^n/G_n\equiv S^n(X/G)$.)

The Macdonald type equation for the orbifold Euler characteristics
of higher orders was obtained in~\cite{Tamanoi}:
\begin{equation}\label{eqn:tamanoi}
 1+\sum_{n=1}^{\infty}\chi^{(k)}(X^n, G_n)\cdot t^n=
 \left(
 \prod_{r_1\ldots,r_k\ge 1}
(1-t^{r_1r_2\cdots r_k})^{r_2r_3^2 \cdots r_k^{k-1}}
\right)^{-\chi^{(k)}(X,G)}.
\end{equation}
For $k=1$ (i.\,e., for the usual orbifold Euler characteristic) the equation~(\ref{eqn:tamanoi}) was obtained in~\cite{Wang}. For the trivial group $G$ it was obtained in~\cite{Bryan_Fulman-1998}.

The aim of this paper is to explain the most part of the computations in~\cite{Tamanoi},
and also of their analogues for the $A$-Euler characteristics
can be replaced by simplest statements,
whose meaning essentially consists in the fact that the invariants
under consideration are invariants of orbifolds.
These statement permit, in particular, to reduce the proof of the
Tamanoi equation to the case of actions of (finite) groups on
the one-point space, i.\,e., for the $G$-spaces $G/G$.
We show that an analogue of the Bryan-Fulman equation (that is, the Tamanoi equation~(\ref{eqn:tamanoi}) for the trivial group $G$)
holds for the $A$-Euler characteristic with an arbitrary finitely
generated group $A$, while an analogue of the Tamanoi equation
does not hold in general.
Thus one can say that the Tamanoi equation is much more strong
(more selective) than the Bryan-Fulman statement.

For the definition of an Euler characteristic (even of the usual one)
one has to assume that the spaces under consideration are nice enough. 
In the framework of this paper, we shall assume that all the spaces
under consideration are invariant unions of cells in finite
$G$-CW-complexes (or rather homeomorphic to such units), see, e.\,g.,~\cite{TtD}. 
(For a finite group $G$, a $G$-CW-complex is a CW-complex
with a cell action of the group $G$ such that, if an element $g\in G$
transforms a cell into itself, then its action on this cell is trivial,
i.\,e., all the points of the cell are fixed by $g$.)
In this form one can represent, in particular, quasiprojective varieties (i.\,e., differences of projective varieties) with (algebraic) actions of finite groups.

\section{The induction property for orbifold Euler characteristics}\label{sec:orbifold}
As in Equation~(\ref{eqn:chi-orb}), initially the orbifold Euler
characteristic was defined for a manifold with a finite group action.
The fact that it was called {\em orbifold} (and in the paper~\cite{HH}
even the Euler number of an orbifoid), means, of course, that the
authors understood that it was an invariant not only of $G$-spaces, but of orbifolds.
However, as far as I understand, a formal definition for orbifolds themselves was first given in~\cite{Roan-1996}.
(The fact that for an orbifold being a full quotient $X/G$, this
definition gives $\chiorb(X,G)$ required a separate proof which was
given in~\cite[Lemma~1]{Roan-1996}.)
However, when we needed the statement (Theorem~\ref{theo:Euler_induction} below), which essentially underlies this fact (or, more precisely, its analogue in a slightly more general context) in the~\cite{GLM-PEMS-2019}, we did not find in the literature not only its proof, but even its formulation.
Because of that we reproduced it there in Teorem~1 (for generalized
(``motivic'') orbifold characteristics of higher orders).

Let us formulate this statement (also in a somewhat more wide context than in~\cite{GLM-PEMS-2019}: not only for the orbifold Euler
characteristics of higher orders) and show that (in the classical,
not in the generalyzed case) it follows from more simple arguments than there.

Let $G$ be a finite group, let $Z$ be a $G$-space (for example, a $G$-CW-complex), and let $H$ be a finite group containing $G$. There is 
a natural {\em induction operation} which creates an $H$-space from $Z$.
Let us consider the following eqivalence relation on the product
$H\times Z$: $(h_1, x_1)\sim(h_2, x_2)$ ($x_i\in Z$, $h_i\in H$) if and only if there exists an element $g\in G$ such that 
$h_2=h_1g^{-1}$, $x_2=gx_1$. On the quotient $\ind_G^H Z:=(H\times Z)/\sim$, one has a natural action of the group $H$:
$f*(h,x)=(fh,x)$ $f,\,h\in H$, $x\in Z$. As a topological space
$\ind_G^H Z$ is homeomorphic to $Z\times H/G$.
The map $x\mapsto (1,x)$ gives an embedding of $Z$ into $\ind_G^H Z$
(considered as a $G$-space).
If $G\subset H\subset K$, one has
$\ind_G^K=\ind_H^K\circ\ind_G^H$. In particular, one has
$\ind_G^H G/G=H/G$, $\ind_H^K H/G=K/G$.

Let $A$ be a finitely generated group.
In order the orbifold $A$-Euler characteristic
(in particular, the usual orbifold Euler characteristic and its analogues of higher orders) to be an invariant of an orbifold
(obviously an additive one),
one needs that, for different partitions of an orbifold into pieces
such that each of them is covered by one uniformizing system
(and therefore is the quotient of a manifold by a finite group action),
the sum of the $A$-Euler characteristics for this pieces does not
depend on the partition. It is not difficult to see that this fact
is equivalent to the one that, for a subgroup $G$ of a (finite)
group $H$ and a $G$-manifold $Z$ one has the equality
$\chi^{(A)}(\ind_G^H Z,H)=\chi^{(A)}(Z,G)$.

\begin{theorem}\label{theo:Euler_induction}
 Let $Z$ be a $G$-invariant union of cells in a finite $G$-CW-complex,
 let $G\subset H$ ($H$ is a finite group), and let $A$ be a finitely
 generated group. Then one has the equality
 \begin{equation}\label{eqn:Euler_induction}
  \chi^{(A)}(\ind_G^H Z,H)=\chi^{(A)}(Z,G)\,.
 \end{equation}
\end{theorem}

\begin{proof} 
 The space $Z$ is the disjoint union of the $G$-orbits of open cells.
 Each of this orbits is isomorphic to
 $\ind_K^G(\sigma^q,K)=\sigma^q\times (G/K)$, where $K$ is a subgroup of $G$,
 $\sigma^q$ is an (open) cell of dimension $q$ with the trivial
 action of the group $K$. Therefore it is sufficient to prove
 Equation~(\ref{eqn:Euler_induction}) for $Z=\sigma^q\times (G/K)$
 (with the trivial action of the group $K$ on $\sigma^q$). Since in this case one has $\ind_G^H Z=\sigma^q\times (H/K)$,
 it is sufficient to show that
 \begin{equation}\label{eqn:for_cell}
 \chi^{(A)}(\sigma^q\times (H/K),H)=
 \chi^{(A)}(\sigma^q\times (K/K),K)\,.
 \end{equation}
 Since the cell $\sigma^q$ of dimension $q$ can be represented as 
 the union of two cells of dimension $q$ and one cell of dimension
 $q-1$, the additivity of $\chi^{(A)}$ implies that
 $\chi^{(A)}(\sigma^q\times (H/K),H)=2\cdot \chi^{(A)}(\sigma^q\times (H/K),H)+\chi^{(A)}(\sigma^{q-1}\times (H/K),H)$ and therefore 
 $\chi^{(A)}(\sigma^q\times (H/K),H)=
 -\chi^{(A)}(\sigma^{q-1}\times (H/K),H)$.
 Therefore it is sufficient to prove~(\ref{eqn:for_cell})
 for $q=0$: 
 \begin{equation}\label{eqn:for_cell0}
 \chi^{(A)}(H/K,H)=
 \chi^{(A)}(K/K,K)\,.
  \end{equation}
 (One could simply use Equation~(\ref{eqn:multiplicat}):
 $\chi^{(A)}(C\times H/K,\{e\}\times H)=
 (-1)^q \cdot\chi^{(A)}(H/K,H)$, where $(-1)^q-\chi(\sigma^q)=
 \chi^{(\{e\})}(\sigma^q,\{e\})$. We have carried out the reduction
 described here in order not to repeat in detail its analogue
 in the proof of Teorem~\ref{theo:Tamanoi_through_pt}.)
 
 One has $\chi^{(A)}(K/K,K)=
 \frac{1}{\vert K\vert}\vert {\rm Hom}(A,K)\vert$.
 Let
 $$
 W=\{(\varphi, [h])\in {\rm Hom}(A,H)\times (H/K): [h]\in (H/K)^{\varphi(A)}\}.$$ and let
 $\pi_1:W\to {\rm Hom}(A,H)$ and
 $\pi_2:W\to H/K$ be the projections to the corresponding factors.
 One has
 \begin{equation}\label{eqn:W_pi_1}
 \vert W\vert=
 \sum_{\varphi\in{\rm Hom}(A,H)}
 \vert \pi_1^{-1}(\varphi)\vert=
 \sum_{\varphi\in{\rm Hom}(A,H)}
 \vert (H/K)^{\varphi(A)}\vert
 =\vert H\vert\cdot\chi^{(A)}(H/K, H).
 \end{equation}
 The isotropy group $H_{[e]}$ of the class of the unit element in $H/K$
 coincides with the subgroup $K$. 
 For $[h]\in H/K$, one has 
 $H_{[h]}=hKh^{-1}$. Hence
 $\pi_2^{-1}([h])={\rm Hom}(A,hKh^{-1})
 \cong {\rm Hom}(A,K)$. Therefore
 \begin{eqnarray}
 \vert W\vert&=&
 \sum_{[h]\in H/K}
 \vert \pi_2^{-1}([h])\vert=
 \sum_{[h]\in H/K}
 \vert {\rm Hom}(A,hKh^{-1})\vert=
 \nonumber\\
 &=&\vert H/K\vert\cdot 
 \vert {\rm Hom}(A,K)\vert
 =\vert H\vert
 \cdot \chi^{(A)}(K/K,K)\,.\label{eqn:W_pi_2}
 \end{eqnarray}
 Comparing (\ref{eqn:W_pi_1}) with (\ref{eqn:W_pi_2}), we get
 (\ref{eqn:for_cell0}). 
\end{proof}

The Theorem means, in particular, that $\chi^{(A)}(\cdot, \cdot)$
is well-defined as a function on the Grothendieck ring
$K_0^{\rm fGr}(\Var)$ of complex quasiprojective varieties
with actions of finite groups: \cite{GLM-PEMS-2019}.

\section{An analogue of the Kapranov zeta function}\label{sec:Tamanoi_eqn}
For a (nice enough) $G$-space $X$, let us denote by
$\zeta^{(A)}_{(X,G)}(t)$ the power series
$$
1+\sum_{n=1}^{\infty} \chi^{(A)}(X^n,G_n)\cdot t^n\in 1+t\cdot
\Q[[t]].
$$
(If $X$ is a complex quasiprojective variety, one has
$$
\zeta^{(A)}_{(X,G)}(t)=\chi^{(A)}\left(\zeta_{(X,G)}(t)\right)\,,
$$
where $\zeta_{(X,G)}(t)\in 1+t\cdot K_0^{\rm fGr}(\Var)[[t]]$ is
the Kapranov zeta function of the pair $(Z, G)$: \cite{GLM-PEMS-2019}.)
Let us denote
$\zeta^{(\Z^{k+1})}_{(X,G)}(t)$ by
$\zeta^{k)}_{(X,G)}(t)$.

In these terms, in~\cite{Bryan_Fulman-1998} it was shown that the series
$\zeta^{(k)}_{(X,\{e\})}(t)$ is equal to a certain series not dependent
on the space (in fact~---
$\zeta^{(k)}_{({\rm pt},\{e\})}(t)$) in the exponent
$\chi(X)\,\left(=\chi^{(k)}(X,\{e\})\right)$.
Besides that the paper contains a computation of the series
$\zeta^{(k)}_{({\rm pt},\{e\})}(t)$. Let us show that the first
statement and even its analogue for the $A$-Euler characteristic
$\chi^{(A)}(\cdot,\cdot)$ with an arbitrary finitely generated group $A$
is a simple consequence of the property~(\ref{eqn:Euler_induction}) (Theorem~\ref{theo:Euler_induction}).
For the proof we shall use the following statement 
necessary in what follows as well.

\begin{lemma}\label{lemma:product}
 Let $X$ and $Y$ be $G$-invariant unions of cells in a finite
 $G$-CW-complex, $Y\subset X$. Then
 \begin{equation}\label{eqn:multiplicativity_zeta}
   \zeta^{(A)}_{(X,G)}(t)=
   \zeta^{(A)}_{(Y,G)}(t)\cdot
   \zeta^{(A)}_{(X\setminus Y,G)}(t)\,.
  \end{equation}
\end{lemma}

 \begin{proof}
 Let $I_0=\{1,\cdots, n\}$, for $I\subset I_0$ by $\overline{I}$ one denote the complement
 $I_0\setminus I$. One has
  $$
  X^n=\bigsqcup_{I\subset I_0}
  \left(Y^I\times (X\setminus Y)^{\overline{I}}\right)=
  \bigsqcup_{s=0}^n\bigsqcup_{\stackrel{I\subset I_0:}{\vert I\vert=s}}
  \left(Y^I\times (X\setminus Y)^{\overline{I}}\right),
  $$
  where the summand
  $$
  \bigsqcup_{\stackrel{I\subset I_0:}{\vert I\vert=s}}
  \left(Y^I\times (X\setminus Y)^{\overline{I}}\right)
  $$
  is a $G_n$-space. Moreover
  $$
  \bigsqcup_{\stackrel{I\subset I_0:}{\vert I\vert=s}}
  \left(Y^I\times (X\setminus Y)^{\overline{I}}\right)\cong
  {\rm Ind}_{G_s\times G_{n-s}}^{G_n}
  \left(Y^{\{1,\ldots, s\}}\times (X\setminus Y)^{\{s+1,\ldots, n\}}\right),
  $$
  where $Y^{\{1,\ldots, s\}}\times (X\setminus Y)^{\{s+1,\ldots, n\}}$ is endowed with a natural $(G_s\times G_{n-s})$-action.
  Therefore
  \begin{eqnarray*}  
  \chi^{(A)}(X^n,G_n)&=&\sum_{s=0}^n
  \chi^{(A)}\left(\bigsqcup_{\stackrel{I\subset I_0:}{\vert I\vert=s}}
  \left(Y^I\times (X\setminus Y)^{\overline{I}}\right),G_n\right)=\\
  &=&\sum_{s=0}^n \chi^{(A)}(Y^{\{1,\ldots, s\}}\times (X\setminus Y)^{\{s+1,\ldots, n\}}, G_s\times G_{n-s})=\\
  &=&\sum_{s=0}^n 
  \chi^{(A)}(Y^s,G_s)\cdot \chi^{(A)}((X\setminus Y)^{n-s}, G_{n-s}),
  \end{eqnarray*}
  what implies (\ref{eqn:multiplicativity_zeta}).
 \end{proof}
 
  \begin{proposition}\label{prop:BrF_for_(A)}
   Let $A$ be a finitely generated and let $Z$ be a union of cells in
   a finite CW-complex. Then one has the equation
 \begin{equation}\label{eqn:BrF_for(A)}
  \zeta^{(A)}_{(Z,\{e\})}(t)=\left(\zeta^{(A)}_{({\rm pt},\{e\})}(t)\right)^{\chi(Z)}.
 \end{equation}
 \end{proposition}
 
 \begin{proof}
  Since $Z$ is a union of cells, due to Lemma~\ref{lemma:product}
  it is sufficient to prove the statement for one cell:
  $$
  \zeta^{(A)}_{(\sigma^q,\{e\})}(t)=\left(\zeta^{(A)}_{({\rm pt},\{e\})}(t)\right)^{(-1)^q}.
  $$
  For $q=1$ the statement is tautological. Since the cell of
  dimension $q$ can be represented as the union of two cells of
  dimension $q$ and one cell of dimension $q-1$, one has
  $$
 \zeta^{(A)}_{(\sigma^q,\{e\})}(t)=\left(\zeta^{(A)}_{(\sigma^{q-1},\{e\})}(t)\right)^{-1}.
 $$
 \end{proof}

\begin{definition}
We shall say that, for a (finitely generated) group $A$, one has
the Tamanoi type equation if for any (nice enough in the sense described
above) $G$-space $X$ one has the equation
\begin{equation}\label{eqn:Tamanoi_for_A}
\zeta^{(A)}_{(X,G)}(t)=
\left(Z_A(t)\right)^{\chi^{(A)}(X,G)},
\end{equation}
where the power series $Z_A(t)$ depends only on the group $A$.
\end{definition}

\begin{remark}
 It is clear that in this case
 $Z_A(t)=\zeta^{(A)}_{({\rm pt},\{e\})}(t)$.
\end{remark}

\begin{theorem}\label{theo:Tamanoi_through_pt}
 For a group $A$, the Tamanoi type equation holds if and only if
 Equation~(\ref{eqn:Tamanoi_for_A}) holds for any (finite) group $K$
 acting (in the unique possible way) on the one-point space ${\rm pt}$, i.\,e.\ if
 \begin{equation}\label{eqn:Tamanoi_for_A_pt}
\zeta^{(A)}_{(K/K,K)}(t)=
\left(Z_A(t)\right)^{\chi^{(A)}(K/K,K)},
\end{equation}
\end{theorem}

\begin{proof}
 Assume that 
 $$
 \zeta^{(A)}_{(K/K,K)}(t)=
\left(Z_A(t)\right)^{\chi^{(A)}(K/K,K)}
 $$
 for any finite group $K$. A $G$-space $X$ (a $G$-invariant union
 of cells in a finite $G$-CW-complex) is the (disjoint) union of
 the orbits $\sigma^q\times(G/K)$ of (open) cells $\sigma^q$. 
 Lemma~\ref{lemma:product} implies that, if $\zeta^{(A)}_{(Y,G)}(t)=
\left(Z_A(t)\right)^{\chi^{(A)}(Y,G)}$ and
$\zeta^{(A)}_{(X\setminus Y,G)}(t)=
\left(Z_A(t)\right)^{\chi^{(A)}(X\setminus Y,G)}$, then $\zeta^{(A)}_{(X,G)}(t)=
\left(Z_A(t)\right)^{\chi^{(A)}(X,G)}$.
 Therefore, it is sufficient to prove Equation~(\ref{eqn:Tamanoi_for_A})
 for the $G$-spaces
 $(\sigma^q\times(G/K), G)$, where $\sigma^q$ is the $q$-dimensional
 cell with the action of the trivial group.
 The reduction used above (based on the fact that the cell of
 dimension $q$ can be represented as the union of two cells of dimension
 $q$ and one cell of dimension $q-1$)
 alongside with Lemma~\ref{lemma:product}
 permit to reduce the proof to the case of the $G$-space $((G/K), G)$.
 Since
 $\left((G/K)^n,G_n\right)=
 {\rm ind}_{K_n}^{G_n}\left((K/K)^n,K_n\right)$,
 one has
 \begin{eqnarray*}
  \zeta_{(G/K,G)}(t)&=&1+\sum_{n=1}^{\infty}
  \chi^{(A)}\left((G/K)^n,G_n\right)\cdot t^n=
  1+\sum_{n=1}^{\infty}
  \chi^{(A)}\left((K/K)^n,K_n\right)\cdot t^n=\\
  {\ }&=&\zeta_{(K/K,K)}(t)=
  (Z_A(t))^{\chi^{(A)}(K/K,K)}=
  (Z_A(t))^{\chi^{(A)}(G/K,G)},
 \end{eqnarray*}
 Q.E.D.
 \end{proof}

\section{Proof of the Tamanoi equation}\label{sec:Tamanoi_proof}
According to Theorem~\ref{theo:Tamanoi_through_pt}, it is sufficient to prove Equation~(\ref{eqn:tamanoi}) for one-point spaces.

\begin{theorem}\label{theo:tamanoi_proof}
 For a finite group $G$ one has
 \begin{equation}\label{eqn:tamanoi_proof}
 \zeta_{(G/G, G)}(t)=
 \left(\prod_{r_1\ldots,r_k\ge 1}
(1-t^{r_1\cdot r_2\cdot\ldots\cdot r_k})^{r_2\cdot r_3^2 \cdot\ldots\cdot r_k^{k-1}}
\right)^{-\chi^{(k)}(G/G,G)},
\end{equation}
(where $\chi^{(k)}(G_n/G_n,G_n)=\frac{1}{\vert G_n\vert}\vert {\rm Hom}(\Z^{k+1}, G_n)\vert$, $\chi^{(k)}(G/G,G)=\frac{1}{\vert G\vert}\vert {\rm Hom}(\Z^{k+1}, G)\vert$).
\end{theorem}

\begin{proof}
 The proof, of course will follow the line of the proof of 
 Equation~(\ref{eqn:tamanoi}) in~\cite{Tamanoi}. 
 The simplification will consist in the fact that we consider
 one-point spaces and therefore one does not need to take care
 of the fixed point sets of elements and of their Euler characteristics.

 We shall use a particular case of Lemma~4-1 from~\cite{Tamanoi}
 which, in the case under consideration is almost obvious.
 Let $K$ be a finite group and let a group $K'$ be an extention
 of $K$ by an element $a$ commuting with all the elements of $K$ and
 such that $a^r\in K$, 
 $\langle a\rangle \cap K=\langle a^r\rangle $ ($r\in \Z_{>0}$).
 
\begin{lemma}\label{lemma3} 
One has the equation
$$
\chi^{(k)}(K'/K', K')=r^k\cdot \chi^{(k)}(K/K, K)\,.
$$
\end{lemma}

\begin{proof}
One has $\chi^{(k)}(K/K, K)=\frac{1}{\vert K \vert}\vert {\rm Hom\,}(\Z^{k+1},K)\vert$, $\chi^{(k)}(K'/K', K')=\frac{1}{\vert K' \vert}\vert {\rm Hom\,}(\Z^{k+1},K')\vert$, $\vert K'\vert=r\cdot\vert K'\vert$. It is obvious that an arbitrary homomorphism $\Phi$ from $\Z^{k+1}$ to
$K'$ is given by the formula $\Phi(e_i)=a^{s_i}\phi(e_i)$, where $e_i$, $i=0,1,\ldots, k+1$, are (free) generators of the group $\Z^{k+1}$,
$\phi$ is a homomorphism from $\Z^{k+1}$ to $K$, $0\le s_i \le r-1$.
Therefore $\vert {\rm Hom\,}(\Z^{k+1},K')\vert=r^{k+1}\cdot
\vert {\rm Hom\,}(\Z^{k+1},K)\vert$, what proves the statement.
\end{proof}
 
 For the proof of~(\ref{eqn:tamanoi}) we shall use the induction on
 the order $k$ of the orbifoid Euler characteristic. For $k=0$ this
 is the usual Macdonald equation. Assune that~(\ref{eqn:tamanoi})
 is proved for the orbifold Euler characteristic of order $k-1$.
 According to Proposition~\ref{prop:A-Euler_for_product} one has
 \begin{equation}\label{eqn;start}
 \zeta_{(G/G, G)}(t)=
 \sum_{n=0}^{\infty}\chi^{(k)}({\rm pt}, G_n)\cdot t^n=\sum_{n=0}^{\infty}\left(
 \sum_{[a]\in {\rm Conj\,}G_n}
 \chi^{(k-1)}({\rm pt}, C_{G_n}(a))\right)\cdot t^n\,.
 \end{equation}
 Conjugacy classes of elements of the group $G_n$ and their
 centralizers are described, in particular, in~\cite[Proposition~1, Lemma~2]{Wang} and in~\cite[Section~3]{Tamanoi}. For $a=({\mathbf{g}},s)\in G_n$ (${\mathbf{g}}=(g_1,\ldots, g_n)$, $s\in S_n$), let us consider the cycles of the permutation $s$. 
 For a cycle $z=(i_1,\ldots, i_r)$ of $s$ its {\em cycle product} is
 $g_{i_1}\cdot\ldots\cdot g_{i_r}$.
 (The cycle product is well-defined up to conjugation.)
 For $[c]\in{\rm Conj\,}G$, let us denote by $m_r(c)$ the number
 of cycles of length $r$ in the permutation $s$ with the cycle
 product from $[c]$.  (One has 
 $\sum\limits_{[c]\in{\rm Conj\,}G,r\ge 1}rm_r(c)=n$.) The collection
 $\{m_r(c)\}$ is called the type of the element
 $a=({\mathbf{g}},s)\in G_n$. Two elements of the group $G_n$
 are conjugate if and only if their types coincide.
 The centralizer $C_{G_n}(a)$ of the element $a=({\mathbf{g}},s)$ is isomorphic to
 $$
 \prod_{[c]\in{\rm Conj\,}G}\prod_{r\ge 1}\left(C_G(c)\cdot\langle a_{r,c}\rangle\right)_{m_r(c)}\,.
 $$
 where $a_{r,c}$ is an element of the group $G_n$, represented by the pair $(\mathbf{g}_c,z)$, where
 $z$ is the cycle $(i_1,\ldots, i_r)$,
 $\mathbf{g}_c=(h_1,\ldots, h_n)$ with $h_{i_j}=g_{i_j}$, $j=1, \ldots, r$,
 and $h_i=e$ (the unit) for all other $i$. Therefore the
 expression~(\ref{eqn;start}) is equal to
 \begin{eqnarray*}
  &{\ }&\sum_{n\ge 0}
  \left(\sum_{\{m_r(c):\sum rm_r(c)=n\}}
  \prod_{[c],r}\chi^{(n-1)}
  \left({\rm pt},(C_G(c)\cdot\langle a_{r,c}\rangle)_{m_r(c)} \right)\right)\cdot t^n=\\
  &=&
  \sum_{\{m_r(c)\}}\left(
  \prod_{[c],r}\chi^{(n-1)}
  \left({\rm pt},(C_G(c)\cdot\langle a_{r,c}\rangle)_{m_r(c)} \right)\right)\cdot t^{\sum rm_{r,c}}=\\
  &=&\prod_{[c],r}\sum_{m_{r,c}=0}^{\infty}
  \chi^{(n-1)}
  \left({\rm pt},(C_G(c)\cdot\langle a_{r,c}\rangle)_{m_r(c)} \right)\cdot t^r=\\
  &=&\prod_{[c],r}
  \prod_{r_1,\ldots,r_{k-1}\ge 1}\left(
  (1-t^{r\cdot r_1\cdot\ldots\cdot r_{k-1}})^{r_2\cdot r_3^2\cdot\ldots\cdot r_{k-1}^{k-2}}
  \right)^{-\chi^{(k-1)}({\rm pt}), C_G(c)\cdot\langle a_{r,c}\rangle)}=\\
  &{\ }&{\text{(here we used the induction assumption)}}\\
  &=&\left(\prod_{r_1,\ldots,r_{k-1}\ge 1}
  (1-t^{r_1\cdot\ldots\cdot r_{k-1}\cdot r})^{r_2\cdot r_3^2\cdot\ldots\cdot r_{k-1}^{k-2}}
  \right)^{-r^{k-1}\sum\limits_{[c]\in {\rm Conj\,}G}\chi^{(k-1)}({\rm pt}), C_G(c))}=\\
  &{\ }&{\text{(here we used Lemma~\ref{lemma3})}}\\
  &=&\left(
  \prod_{r_1,\ldots,r_{k-1}, r_k\ge 1}
  (1-t^{r_1\cdot\ldots\cdot r_k })^{r_2\cdot r_3^2\cdot\ldots\cdot r_{k}^{k-1}}
  \right)^{-\chi^{(k)}({\rm G/G}, G)}\,.
 \end{eqnarray*}
\end{proof}
 
 Let us show that the Tamanoi equation does not hold, in general,
 for groups $A$ different from free abelian ones.
 One has
$$
\zeta^{(\Z_2)}_{({\rm pt},\{e\})}(t)=1+t+t^2+\cdots,
$$
$$
\zeta^{(\Z_2)}_{({\rm pt},\Z_2)}(t)=1+t+\frac{3}{4}\cdot t^2+\cdots,
$$
$\chi^{(\Z_2)}({\rm pt}, \Z^2)=1$. Thus
$\zeta^{(\Z_2)}_{({\rm pt},\Z_2)}(t)\ne\left(\zeta^{(\Z_2)}_{({\rm pt},\{e\})}(t)\right)^{\chi^{(\Z_2)}({\rm pt}, \Z^2)}$.
This example may look not convincing enaugh since $\chi^{(\Z_2)}(\cdot, \cdot)$ is not an integer valued invariant. In particular, in~\cite{Tamanoi}, the $A$-Euler characteristics
are defined only for groups $A$, with the group $\Z$ as a
direct summand. Let us give an appropriate example. One has
$$
\zeta^{(\Z\times\Z_2)}_{({\rm pt},\{e\})}(t)=1+t+2t^2+ \cdots,
$$
$$
\zeta^{(Z\times\Z_2)}_{({\rm pt},\Z_2)}(t)=1+2t+4t^2+\cdots,
$$
$\chi^{(\Z\times\Z_2)}({\rm pt}, \Z_2)=2$. Therefore
$\zeta^{(\Z\times\Z_2)}_{({\rm pt},\Z_2)}(t)\ne\left(\zeta^{(\Z\times\Z_2)}_{({\rm pt},\{e\})}(t)\right)^{\chi^{(\Z\times\Z_2)}({\rm pt}, \Z^2)}$.

\end{document}